\documentclass[10pt,draft,reqno]{amsart}
     \makeatletter
     \def\section{\@startsection{section}{1}%
     \z@{.7\linespacing\@plus\linespacing}{.5\linespacing}%
     {\bfseries%\normalfont\scshape
     \centering
     }}
     \def\@secnumfont{\bfseries}
     \makeatother
\newtheorem{theorem}{Theorem}[section]
\newtheorem{lemma}[theorem]{Lemma}

\newtheorem{corollary}[theorem]{Corollary}
\theoremstyle{definition}
\newtheorem{definition}[theorem]{Definition}

\theoremstyle{remark}

\numberwithin{equation}{section} 
 \newcommand{\po}{\partial\Omega}
 
\begin{document}

\title[]
 {Dependence Result of the Weak Solution of Robin Boundary Value Problems }

\author{Akhlil Khalid}

\address{
Insitut of Applied Analysis
\\
Ulm University
\\
89069 Ulm
\\
Germany}

\email{akhlil.khalid@gmail.com}

\thanks{This work was partially financed by DAAD (Deutscher Akademischer Austausch Dienst)}

\thanks{}



\date{}

\dedicatory{}

\commby{}

%%% ----------------------------------------------------------------------

\begin{abstract}
In this article we establish an approximation result involving the Laplacian with Robin boundary conditions. It 
informs about the weak solution's dependence from the input function on the boundary.

\end{abstract}

%%% ----------------------------------------------------------------------
\maketitle
%%% ----------------------------------------------------------------------
\section{Introduction}

Let $\Omega$ be a bounded domain with Lipschitz boundary. We consider the problem of the Laplacian with Robin
boundary conditions,
\begin{equation}
 \frac{\partial u}{\partial \nu}+\beta u=0
\end{equation}
where $\nu$ is the outwar normal verctor and $\beta$ is a measurable positive bounded function on the boundary
$\po$. This kind of problems was extensively studied by many autors, we refer to \cite{A1}, \cite{A2}, \cite{AW1} 
, \cite{War}, \cite{D} and references therein for more details.

The aim of this article is to show a dependance result of a sequence of weak solutions $(u_n)_{n\geq 0}$ 
with a sequence of ´´ input´´ functions $(\beta)_{n\geq 0}$. The proof is based on a technical Lemma due to
Stampaccia \cite{Sta}.
\section{Preliminaries and main result}

We assume that $\Omega\subset\mathbb R^d$ $(d\geq 3)$ is a bounded domain with Lipschitz boundary. We denote 
by $\sigma$ the restriction to $\partial\Omega$ of the $(d-1)-$dimentional Hausdorff measure.

We know that the following continous embedding holds,

\begin{equation}\label{emb1}
 H^1(\Omega)\rightarrow L^{q}(\Omega),\text{  } q=\dfrac{2d}{d-2}
\end{equation}

Moreover each function $u\in H^1(\Omega)$ has a trace which is in $L^s(\partial\Omega)$, where $s=\frac{2(d-1)}{d-2}$; i.e. there is a constant $c>0$ such that

\begin{equation}\label{emb2}
 \|u\|_{s,\partial\Omega}\leq c \|u\|_{H^1(\Omega)}\text{ for all } u\in H^1(\Omega)
\end{equation}

Let $\lambda>0$ be a real number, $f\in L^p(\Omega)$ $(p>d)$ and $\beta$ be a nonnegative bounded measurable
function on $\po$. We consider the following Robin boundary value problem

\begin{equation}\label{eq:eq}
 \begin{cases}
  -\Delta u+\lambda u=f& \text{ in }\Omega\\[0.2cm]
  \frac{\partial u}{\partial\nu}+\beta u=0& \text{ in }\po
 \end{cases}
\end{equation}

The form associated with the Laplacian with Robin boundary condition is
\[
a_{\beta}(u,v)=\int_{\Omega}\nabla u\nabla v dx+ \int_{\partial\Omega}
\beta u vd\sigma \text { for all } u,v\in H^1(\Omega)
\]

We start by the definition of the weak solution of the problem \eqref{eq:eq}.
\begin{definition}
Let $f\in L^p(\Omega)$. For each $\lambda>0$, a function $u=G_{\beta}^{\lambda}f\in H^1(\Omega)$ is
called a weak solution of the Robin boundary value Problem (associated with $\beta$) if for every $v\in H^1(\Omega)$

\[
 a^{\lambda}_{\beta}(u,v)=\int_{\Omega}fvdx,
\]
where for $u,v\in H^1(\Omega)$ 
\[
 a_{\beta}^{\lambda}(u,v)=a_{\beta}(u,v)+\lambda\int_{\Omega}uv\ dx
\]
\end{definition}

It is clear that the closed bilinear form $a_{\beta}$ is continous on $H^1(\Omega)$ and also coercive on 
$H^1(\Omega)$ in the sens that there exists a constant $c>0$ such that for all $u\in H^1(\Omega)$
\[
 a^{\lambda}_{\beta}(u,u)\geq \|u\|^2_{H^1(\Omega)}
\]

Let $L$ be the linear functional on $H^1(\Omega)$ defined by : for $v\in H^1(\Omega)$
\[
 Lv:=\int_{\Omega}fv\ dx
\]

Since $p\geq 2$, the functional $L$ is well defined and continous on $H^1(\Omega)$. Thus by coerciveness of the 
bilinear form $a_{\beta}$, the Lax-Milgram Lemma (see \cite[Corollaire V.8 p:84]{Br}) implies that there exists a
unique weak solution $u\in H^1(\Omega)$ of the boundary value probem \eqref{eq:eq}.

The following lemma is important in the proof of Theorem \ref{thm1}, we can find its proof in \cite{Sta} Lemma 4.1.

\begin{lemma}\label{lem1}

Let $\varphi=\varphi(t)$ be a nonnegative, nonincreasing function on the half line $t\geq k_0\geq 0$ such that there are positive constants $c,\alpha$ and $\delta(\delta>0)$ such that 
\[
\varphi(h)\leq c (h-k)^{-\alpha}\varphi(k)^{\delta}
\]
for all $h>k\geq k_0$. Then we have
\[
\varphi(k_0+d)=0\text{, where } d>0 \text{ satisfies } d^{\alpha}=c\varphi(k_0)^{\delta-1}2^{\delta(\delta-1)}
\]

 \end{lemma}

\begin{theorem}\label{thm0}
Let $u$ be a weak solution and assume that $p>d$.
Then 

1) if $\lambda=0$ and $\Omega$ is of finite volume, there exists a strictly positive constants
$C_1=C_1(d,p,|\Omega|)$ such that 
\[
 |u(x)|\leq C_1\|f\|_p\quad\text{    a.e on }\overline{\Omega}
\]

2) if $\lambda>0$ and $\Omega$ is an arbitrary domain, there exist a strictly positive constant $C_2=C_2(d,p,\lambda)$
 such that
\[
 |G^{\lambda}_{\beta}f(x)|\leq C_2\|f\|_p\quad\text{    a.e on }\overline{\Omega}
\]

\end{theorem}

The proof can be found in \cite{War} and is based on the Maza'ya inequality and a standard
argument as in Theorem 4.1 of \cite{Sta}.
\\

Our main result is the following Theorem,
\begin{theorem}\label{thm1}

Any sequence $(u_n)_{n\geq 0}$ of weak solutions of the Robin boundary value problem associated to the sequence $(\beta_n)_{\geq 0}$ verify the following inequality:
\begin{equation}
\|u_n-u_m\|_{\infty,\overline{\Omega}}\leq C\|u_n\|_{\infty,\partial\Omega}\|\beta_n-\beta_m\|_{\infty,\partial\Omega}
\end{equation}
for all $n,m\in \mathbb N$ and where $C$ may depend of $\lambda$.
\end{theorem}

\section{Proof of Theorem \ref{thm1}}
\begin{proof}
 Let $(u_n)_{n\geq 0}$ be a sequence of weak solutions associated with the sequence $(\beta_n)_{\geq 0}$. Let $k\geq 0$ be a real number and define $u_{n,m}:=u_n-u_m$.

Define $v_{n,m}:=(|u_{n,m}|-k)^+\mathrm{sgn}(u_{n,m})$. Then $v_{n,m}\in H^1(\Omega)$ and 
\[
 \nabla v_{n,m} =\left\{
          \begin{array}{ll}
           \nabla u_{n,m} &  \text{in} \quad A_{n,m}(k); \\
            0 &
            \text{otherwise}
          \end{array}
          \right.
\]
where $A_{n,m}(k)=\{x\in\overline{\Omega}: |u_{n,m}(x)|>k\}$. In the following, we write $u,v,A(k)..$ instead of $
u_{n,m},v_{n,m},A_{n,m}(k)..$.

It is clear that $a^{\lambda}_{\beta_n}(u_n,v)-a^{\lambda}_{\beta_m}(u_m,v)=0$. Calculating we obtain:
\begin{equation}
\begin{split}
           0 & = \int_{\Omega}\nabla(u_n-u_m)\nabla v dx+\int_{\partial\Omega}(\beta_nu_n-\beta_mu_m)vd\sigma+\lambda\int_{\Omega}(u_n-u_m)vdx\\
             & = \int_{A(k)}|\nabla v|^2dx+\int_{\partial\Omega}(\beta_n-\beta_m)u_n+\beta_m(u_n-u_m)v d\sigma+\lambda\int_{\Omega}(u_n-u_m)vdx\\
             & = \int_{A(k)}|\nabla v|^2dx+\int_{\partial\Omega\cap A(k)}(\beta_n-\beta_m)u_nvd\sigma+\int_{\partial\Omega\cap A(k)}\beta_m(u_n-u_m)vd\sigma\\
             &+\lambda\int_{A(k)}(u_n-u_m)vdx\\
             & = \int_{A(k)}|\nabla v|^2dx+\int_{\partial\Omega\cap A(k)}(\beta_n-\beta_m)u_nvd\sigma+\int_{\partial\Omega\cap A(k)}\beta_m v^2d\sigma\\
              &+k\int_{\partial\Omega\cap A(k)}\beta_m |v|d\sigma+\lambda\int_{A(k)}v^2dx+\lambda k\int_{A(k)}|v|dx\\
             & = a^{\lambda}_{\beta_m}(v,v)+\int_{\partial\Omega\cap A(k)}(\beta_n-\beta_m)u_nvd\sigma+k\int_{\partial\Omega\cap A(k)}\beta_m |v|d\sigma+\lambda k\int_{A(k)}|v|dx
\end{split}
\end{equation}
It follows that 
\begin{equation}
\begin{split}
 a^{\lambda}_{\beta_m}(v,v)+\int_{\partial\Omega\cap A(k)}(\beta_n-\beta_m)u_nvd\sigma&=-k\int_{\partial\Omega\cap A(k)}\beta_m |v|d\sigma-\lambda k\int_{A(k)}|v|dx\\
&\leq 0
\end{split}
\end{equation}
Which leads to
\[
 a^{\lambda}_{\beta_m}(v,v)  \leq \int_{\partial\Omega\cap A(k)}(\beta_m-\beta_n)u_nvd\sigma
\]
Using the H\"{o}lder inequality and \eqref{emb2}, we obtain the following estimates,
\begin{equation}
\begin{split}
a^{\lambda}_{\beta_m}(v,v) & \leq \int_{\partial\Omega\cap A(k)}(\beta_m-\beta_n)u_nvd\sigma\\
                 & \leq \|\beta_n-\beta_m\|_{\infty,\partial\Omega}\int_{\partial\Omega\cap A(k)}u_nvd\sigma\\
                 & \leq \|\beta_n-\beta_m\|_{\infty,\partial\Omega}\|u_n\|_{2,\partial\Omega\cap A(k)}
                 \|v\|_{2,\partial\Omega\cap A(k)}\\
                 & \leq \|\beta_n-\beta_m\|_{\infty,\partial\Omega}
                 \|u_n\|_{\infty,\partial\Omega}|\partial\Omega\cap A(k)|^{\frac{1}{2}}
                 |\partial\Omega\cap A(k)|^{\frac{1}{2}-\frac{1}{s}}\|v\|_{s,\partial\Omega}\\
                 &  \leq \|\beta_n-\beta_m\|_{\infty,\partial\Omega}
                 \|u_n\|_{\infty,\partial\Omega}|\partial\Omega\cap A(k)|^{1-\frac{1}{s}}\|v\|_{s,\partial\Omega}\\
                 & \leq c\|\beta_n-\beta_m\|_{\infty,\partial\Omega}
                 \|u_n\|_{\infty,\partial\Omega}|\partial\Omega\cap A(k)|^{1-\frac{1}{s}}\|v\|_{H^1(\Omega)}
\end{split}
\end{equation}

We have then,
\begin{equation}
\begin{split}
\alpha \|v\|^2_{H^1(\Omega)} & \leq a^{\lambda}_{\beta_m}(v,v) \\
                           &\leq c\|\beta_n-\beta_m\|_{\infty,\partial\Omega}\|u_n\|_{\infty,\partial\Omega}
                           |\partial\Omega\cap A(k)|^{1-\frac{1}{s}}\|v\|_{H^1(\Omega)}
\end{split}
\end{equation}

It follows that  
\begin{equation}
\begin{split}
 \|v\|_{H^1(\Omega)}&\leq c_1\|\beta_n-\beta_m\|_{\infty,\partial\Omega}
 \|u_n\|_{\infty,\partial\Omega}|\partial\Omega\cap A(k)|^{1-\frac{1}{s}}
\end{split}
\end{equation}

Using the inequalities \eqref{emb1} and \eqref{emb2}, we obtain the following estimates,
\begin{equation}\label{eq1}
\begin{split}
\|v\|_{s,\partial\Omega\cap A(k)} &\leq c_2\|\beta_n-\beta_m\|_{\infty,\partial\Omega}
\|u_n\|_{\infty,\partial\Omega}|\partial\Omega\cap A(k)|^{1-\frac{1}{s}}
\end{split}
\end{equation}
and,
\begin{equation}\label{eq2}
\begin{split}
\|v\|_{q,A(k)} &\leq c_3\|\beta_n-\beta_m\|_{\infty,\partial\Omega}
\|u_n\|_{\infty,\partial\Omega}|\partial\Omega\cap A(k)|^{1-\frac{1}{s}}
\end{split}
\end{equation}

Let now $h>k\geq 0$. Then $A(h)\subset A(k)$ and on $A(h)$ we have $|v|\geq h-k$. It follows that

\begin{equation}\label{eq3}
\begin{split}
\|v\|_{s,\partial\Omega\cap A(k)} & \geq \|v\|_{s,\partial\Omega\cap A(h)}\\
                                  & \geq \||u|-k\|_{s,\partial\Omega\cap A(h)}\\
                                  & \geq (h-k)|\partial\Omega\cap A(h)|^{\frac{1}{s}}
\end{split}
\end{equation}

We deduce from \eqref{eq1} that 
\[
(h-k)|\partial\Omega\cap A(h)|^{\frac{1}{s}}\leq c_2\|\beta_n-\beta_m\|_{\infty,\partial\Omega}
\|u_n\|_{\infty,\partial\Omega}|\partial\Omega\cap A(k)|^{1-\frac{1}{s}}
\]

which reduces to,
\[
|\partial\Omega\cap A(h)|\leq c_2^s(h-k)^{-s}\|\beta_n-\beta_m\|^s_{\infty,\partial\Omega}
\|u_n\|^s_{\infty,\partial\Omega}|\partial\Omega\cap A(k)|^{s-1}
\]

Set $\phi(h)=|\partial\Omega\cap A(h)|$, we obtain,
\[
\phi(h)\leq C (h-k)^{-s}\phi(k)^{s-1}
\]
where $C=c_2^s\|\beta_n-\beta_m\|^s_{\infty,\partial\Omega}\|u_n\|^s_{\infty,\partial\Omega}$.

As $s-1>1$, then the conditions of the Lemma \ref{lem1} are satisfied with $\delta=s-1$ and $k_0=0$, 
one obtain $\phi(d)=0$ where $d>0$ satisfies $d^s=C\phi(0)^{s-2}2^{(s-1)(s-2)}$, consequently 
\[
d=c_4\|\beta_n-\beta_m\|_{\infty,\partial\Omega}\|u_n\|_{\infty,\partial\Omega}                                                 
\]
and 
\begin{equation}\label{eq4}
\|u_n-u_m\|_{\infty,\partial\Omega}\leq c_4\|u_n\|_{\infty,\partial\Omega}
\|\beta_n-\beta_m\|_{\infty,\partial\Omega}
\end{equation}

In the same way as in \eqref{eq3}, we obtain
\[
 \|v\|_{q,A(k)}\geq(h-k)|A(k)|^{\frac{1}{q}}
\]

From \eqref{eq2}, we deduce
\[
(h-k)|A(h)|^{\frac{1}{q}}\leq c_3\|\beta_n-\beta_m\|_{\infty,\partial\Omega}
\|u_n\|_{\infty,\partial\Omega}|\partial\Omega\cap A(k)|^{1-\frac{1}{s}}
\]

We take $k=d$ and $h=\gamma d$ with $\gamma>1$, we obtain $|A(\gamma d)|=0$ which leads to 
\begin{equation}\label{eq5}
\begin{split}
\|u_n-u_m\|_{\infty,\Omega} & \leq \gamma d\\
                           & \leq \gamma c_4\|u_n\|_{\infty,\partial\Omega}
                           \|\beta_n-\beta_m\|_{\infty,\partial\Omega}
\end{split}
\end{equation}

From \eqref{eq4} and \eqref{eq5} we obtain our Theorem.
\end{proof}

\begin{corollary}
 Let $(u_n)_{n\geq 0}$ be a sequence weak solutions associated with the sequence
 $(\beta_n)_{\geq 0}\in L^{\infty}(\partial\Omega)$  such that $\inf_n\beta_n>0$ then if
 $(u_n)_{n\geq 0}$ is uniformly bounded we have 
 for $p>d$
 \begin{equation}
\|u_n-u_m\|_{\infty,\overline{\Omega}}\leq C\|f\|_{p}\|\beta_n-\beta_m\|_{\infty,\partial\Omega}
\end{equation}
for all $n,m\in \mathbb N$ and where $C$ may depend of $\lambda$.

\end{corollary}

In the case where the sequence of weak solutions $(u_n)_{n\geq 0}$ is uniformely bounded with respect to $n$ we have
the following consequence
\begin{corollary}
 Let $(u_n)_{n\geq 0}$ be a sequence weak solutions associated with the sequence
 $(\beta_n)_{\geq 0}\in L^{\infty}(\partial\Omega)$  such that $\inf_n\beta_n>0$  and 
 $\lim_{n}\beta_n(x)=\beta(x)\text{ a.e } x\in\partial\Omega$ then if $(u_n)_{n\geq 0}$ is uniformly bounded we have 
 $\lim_{n}u_n(x)=u(x)\text{ a.e } x\in\overline{\Omega}$, where $u$ is the weak solution associated with $\beta$.
\end{corollary}

%------------------------------------------------------------------

\par\bigskip

%{\bf Acknowledgment.}
%------------------------------------------------------------------------------------------------------------------------

\end{document}